\documentclass{amsart}
\theoremstyle{plain}
\newtheorem{Thm}{Theorem}

\begin{document}
\large
\title[Interior Gradient Bound]{Interior Gradient Bound For Minimal Graphs in a Product Manifold
$M\times R$}

\author{Li MA and Dezhong Chen}

\address{Department of mathematical sciences \\
Tsinghua university \\
Beijing 100084 \\
China}

\email{lma@math.tsinghua.edu.cn} \dedicatory{}
\date{May 26th, 2004}

\thanks{The work is partially supported by the key 973 project of the Ministry of Science and Technology
 of China. Part of the work is done while the first name
author is visiting Greifswald University, Germany. The first name
author thanks Prof. J.Eichhorn for the invitation.}



\begin{abstract}
Let $(M, g)$ be an $n$-dimensional complete Riemannian manifold
with $Ric(M)\geq-(n-1)Q$, where $Q\geq0$ is a constant. We obtain
an interior gradient bound for minimal graphs in $M\times R$ under
some technical assumptions. For details, see Theorem 2.
\end{abstract}
\maketitle

\begin{section}{Introduction}
  Let $(M,g)$ be an $n$-dimensional complete Riemannian manifold.
  For each $u\in C^{3}(M)$, we can define a graph in the product
  manifold $(M\times R,g+dt^{2})$ as follows:
  $$
  \Sigma(u):=\{(x,u(x))|x\in M\}.
  $$
  Naturally, we can equip the graph $\Sigma(u)$ with the  metric induced from
  $M\times R$. Then, it is well known that $\Sigma(u)$ is a minimal
  hyper-surface in the product manifold $M\times R$
  iff $u$ satisfies
  \begin{equation}
  div_{M}\frac{\nabla_{M}u}{\sqrt{1+|\nabla_{M}u|^{2}}}=0,
  \end{equation}
  where $\nabla_{M}$ and $div_{M}$ are the gradient and divergence
  of $(M,g)$. Minimal surface theory in $R^{n+1}$ has been studied
  by many famous mathematicians. Minimal surfaces in a general compact Riemannian
  manifold $(N,\bar{g})$ are also interesting subjects with rich
  applications. There is relatively few results of minimal surfaces
  in a general complete Riemannian manifold except the cases when $N=H^{n+!}(-1)$
  or $M\times R$, see \cite{M1} and  \cite{R}.

  In the case of $M$ being a domain of the Euclidean space $R^{n}$,
  many powerful techniques have
  already been developed to derive a priori gradient bounds for
  solutions to (1). The key tool is the use of the Maximum Principle
  for suitable quantities. In
  \cite{K}, N. Korevaar gave a beautiful application of the maximum principle proof of the
  following
  \begin{Thm}
  Let $M$ be the ball $B_2(0)$ of radius 2 with center at the origin
  $0$  of $R^{n}$. Then there exist two constants $K_{i}$, $i=1,2$, depending
  only on $n$, such that
  $$
  v(0)\leq K_{2}\exp(K_{1}u(0)^{2}),
  $$
  where $u\in C^{3}(B_{1}(0))$ is a negative solution of (1),
  $v=\sqrt{1+|\nabla u|^{2}}$.
  \end{Thm}
Later, this kind of argument was used by Ecker and Huisken
\cite{EH} and Colding and Minicozzi  \cite{CM} to study the mean
curvature flow.

  In this note, we will use the maximum principle argument to prove the
  following
  \begin{Thm}
  Let $M$ be an $n$-dimensional complete Riemannian manifold with
  $Ric(M)\geq-(n-1)Q$, where $Q\geq0$ is a constant. $B_{r}(p)$ is
  a geodesic ball with center $p$ and radius $r$. Assume that
  there exists $0\leq C_{+}<\frac{\pi^{2}}{8r^{2}}$ such that for
  any geodesic $\gamma:[0,s]\rightarrow B_{r}(p)$, $0\leq s<r$,
  parametrized by arc length with $\gamma(0)=p$, we have at
  $\gamma(t)$, $0\leq t\leq s$,
  $$
  K(X,\frac{\partial}{\partial\gamma})\leq C_{+},
  $$
  where $X$ is any vector of length less than one in $T_{\gamma(t)}M$ linearly independent
  of $\frac{\partial}{\partial\gamma}$,
  $K(X,\frac{\partial}{\partial\gamma})$ is the sectional
  curvature of the 2-plane spanned by $X$ and
  $\frac{\partial}{\partial\gamma}$. Then there
  exist $K_{i}=K_{i}(n,Q,r,C_{+})$, $i=1,2$, such that if $u\in
  C^{3}(B_{r}(p))$ is a negative solution of (1), then
  $$
  v(p)\leq K_{2}\exp(K_{1}u(p)^{2}).
  $$
  \end{Thm}

  Our argument is to generalize the proof of Korevaar to our case.
  However, many details are much more involved. We believe the
  estimate in the theorem above can be improved with the use of the idea in \cite{CM}.
  Another interesting subject is to study the Mean Curvature Flow
  in $M\times R$, see \cite{M2}.
\end{section}

\begin{section}{proof of theorem 2}
  Our proof is divided into six steps.

  {\bf Step 1.} Let $\eta(x,z)$ be a nonnegative continuous
  function in $B_{r}(p)\times R^{-}$, which vanishes on
  $\{|x|=r,z<0\}$ and is smooth where it is positive. Then $\eta v$
  has a positive maximum in the interior of $B_{r}$, say at $P$.
  If $\{x^{i}\}$ is a normal coordinate system at $P$ with $\frac{\partial}{\partial
  x^{n}}=\frac{\partial}{\partial \gamma}$, then we have at $P$
  $$
  (\eta v)_{i}=0,i=1,\cdots,n,
  $$
  $$
  [(\eta v)_{ij}]\leq0,
  $$
  where the subscript $i$ on the outside of parenthesis means the
  covariant derivative w.r.t. $\frac{\partial}{\partial x^{i}}$.
  By the chain rule,
  $$
  (\eta)_{i}=\eta_{i}+\eta_{z}u_{i},
  $$
  $$
  (\eta)_{ij}=\eta_{ij}+\eta_{iz}u_{j}+\eta_{zj}u_{i}+\eta_{zz}u_{i}u_{j}+\eta_{z}u_{ij}.
  $$

  Define
  $$
  g^{ij}=\delta^{ij}-\nu^{i}\nu^{j},
  $$
  where $\delta^{ij}$ is the standard Kronecker's symbol and
  $\nu^{i}=\frac{u_{i}}{v}$. Note that $|\nu|<1$. So $[g^{ij}]$ is positive
  definite. Then the trace of the product of $[(\eta v)_{ij}]$
  with $[g^{ij}]$ is non-positive, i.e.,
  $$
  g^{ij}(\eta v)_{ij}\leq 0.
  $$
  But
  $$
  (\eta
  v)_{ij}=(\eta)_{ij}v+(\eta)_{i}v_{j}+(\eta)_{j}v_{i}+\eta v_{ij}.
  $$
  Hence
  $$
  0\geq g^{ij}(\eta)_{ij}v+g^{ij}(\eta)_{i}v_{j}+g^{ij}(\eta)_{j}v_{i}+g^{ij}\eta v_{ij}
  $$
  $$
  \hspace{-9mm}=g^{ij}(\eta)_{ij}v+2g^{ij}(-\frac{\eta}{v})v_{i}v_{j}+\eta g^{ij}v_{ij}
  $$
  $$
  \hspace{-11mm}=v[g^{ij}(\eta)_{ij}+\eta\frac{g^{ij}}{v}(-\frac{2}{v}v_{i}v_{j}+v_{ij})].
  $$
  Here we have substituted $-\frac{\eta}{v}v_{i}$ for
  $(\eta)_{i}$. Therefore
  \begin{equation}
  g^{ij}(\eta)_{ij}+\eta\frac{g^{ij}}{v}(-\frac{2}{v}v_{i}v_{j}+v_{ij})\leq0.
  \end{equation}

  We compute
  $$
  v_{i}=\nu^{j}u_{ji},
  $$
  $$
  v_{ij}=\nu^{k}u_{kij}+\frac{1}{v}(u_{kj}u_{ki}-\nu^{k}u_{ki}\nu^{l}u_{lj}).
  $$
  So
  \begin{equation}
  g^{ij}v_{ij}=g^{ij}u_{kij}\nu^{k}+\frac{1}{v}g^{ij}(u_{kj}u_{ki}-\nu^{k}u_{ki}\nu^{l}u_{lj}).
  \end{equation}
  Note that the second term is nonnegative for $|\nu|<1$. From
  $$
  u_{ki}=u_{ik}
  $$
  and the Ricci formula
  $$
  u_{ikj}=u_{ijk}+R_{kj}u_{i}
  $$
  where $R_{kj}$ is the Ricci curvature of $(M,g)$, it follows
  that
  $$
  g^{ij}u_{kij}\nu^{k}=g^{ij}(u_{ijk}+R_{kj}u_{i})\nu^{k}
  $$
  \begin{equation}
  \hspace{24mm}=g^{ij}u_{ijk}\nu^{k}+\frac{g^{ij}}{v}R_{kj}u_{i}u_{k}.
  \end{equation}
  For $u$ satisfies (1), we have
  $$
  g^{ij}u_{ij}=0.
  $$
  Integrating by parts,
  \begin{equation}
  g^{ij}u_{ijk}=(g^{ij}u_{ij})_{k}-(g^{ij})_{k}u_{ij}=-(g^{ij})_{k}u_{ij}.
  \end{equation}
  But
  $$
  -(g^{ij})_{k}=(\nu^{i}\nu^{j})_{k}=\frac{u_{ik}u_{j}}{v^{2}}+\frac{u_{i}u_{jk}}{v^{2}}
  -\frac{2}{v^{4}}u_{i}u_{j}u_{l}u_{lk}.
  $$
  So
  \begin{equation}
  -(g^{ij})_{k}u_{ij}\nu^{k}=2(\frac{u_{ik}u_{ij}u_{j}u_{k}}{v^{3}}-\frac{u_{i}u_{j}u_{k}u_{l}u_{ij}u_{kl}}{v^{5}})
  =\frac{2}{v}g^{ij}v_{i}v_{j}.
  \end{equation}
  Note that
  \begin{equation}
  g^{ij}R_{kj}u_{i}u_{k}=R_{kj}\nu^{j}\nu^{k}.
  \end{equation}
  Combining (3)-(7) yields
  \begin{equation}
  \frac{1}{v}g^{ij}(-\frac{2}{v}v_{i}v_{j}+v_{ij})\geq\frac{R_{kj}\nu^{j}\nu^{k}}{v^{2}}.
  \end{equation}

  {\bf Step 2.} Let
  $$
  u_{0}=-u(p)>0.
  $$
  Define
  $$
  \phi(x,z)=(\frac{1}{2u_{0}}z+(r^{2}-\rho(x)^{2}))^{+}.
  $$
  Here $"+"$ means positive part and
  $$
  \rho(x)=dist(p,x)
  $$
  for all $x\in B_{r}(p)$, where $dist(\cdot,\cdot)$ denotes the
  geodesic distance. Then
  $$
  0\leq\phi\leq
  r^{2},\phi_{z}=\frac{1}{2u_{0}},\phi_{zz}=\phi_{iz}=0,
  $$
  $$
  |\nabla\phi|^{2}=4\rho^{2}|\nabla\rho|^{2}=4\rho^{2}<4r^{2}.
  $$
  Let
  $$
  f(\phi)=\exp(C_{1}\phi)-1,
  $$
  where $C_{1}>0$ is to be determined. It is easy to see that
  $$
  f(0)=0,f'>0,f''>0.
  $$
  Let
  $$
  \eta(x,u(x))=f\circ\phi(x,u(x)).
  $$
  Then $\eta$ has the properties required at the beginning of {\bf
  Step 1}. In this setting,
  $$
  (\eta)_{i}=f'\cdot(\phi_{i}+\phi_{z}u_{i}),
  $$
  $$
  (\eta)_{ij}=f''\cdot(\phi_{j}+\phi_{z}u_{j})(\phi_{i}+\phi_{z}u_{i})
  $$
  $$
  \hspace{42mm}+f'\cdot(\phi_{ij}+\phi_{iz}u_{j}+\phi_{zj}u_{i}+\phi_{zz}u_{j}u_{i}+\phi_{z}u_{ij}).
  $$
  Hence
  $$
  g^{ij}(\eta)_{ij}=g^{ij}f''\cdot(\phi_{i}\phi_{j}+\phi_{z}\phi_{i}u_{j}+\phi_{z}\phi_{j}u_{i}
  +\phi_{z}^{2}u_{i}u_{j})
  $$
  $$
  \hspace{-19mm}+g^{ij}f'\cdot(\phi_{ij}+\phi_{z}u_{ij})
  $$
  $$
  \hspace{29mm}=f''\cdot(g^{ij}\phi_{i}\phi_{j}+2\phi_{z}g^{ij}\phi_{i}u_{j}+\phi_{z}^{2}g^{ij}u_{i}u_{j})
  +f'\cdot g^{ij}\phi_{ij}.
  $$
  But
  $$
  g^{ij}\phi_{i}\phi_{j}=\phi_{i}^{2}-\phi_{i}\nu^{i}\phi_{j}\nu^{j}\geq0,
  $$
  $$
  g^{ij}\phi_{i}u_{j}=\frac{\phi_{i}\nu^{i}}{v},
  $$
  $$
  g^{ij}u_{i}u_{j}=\frac{|\nabla u|^{2}}{v^{2}}.
  $$
  So
  $$
  g^{ij}(\eta)_{ij}=f''\cdot(\phi_{z}^{2}\frac{|\nabla u|^{2}}{1+|\nabla u|^{2}}+\frac{2\phi_{z}\phi_{i}\nu^{i}}{v}
  +\phi_{i}^{2}-\phi_{i}\nu^{i}\phi_{j}\nu^{j})+f'\cdot g^{ij}\phi_{ij}
  $$
  \begin{equation}
  \hspace{-28mm}\geq f''\cdot\frac{|\nabla u|^{2}+4u_{0}u_{i}\phi_{i}}{4u_{0}(1+|\nabla
  u|^{2})}+f'\cdot g^{ij}\phi_{ij}.
  \end{equation}
  Combining (2), (8) and (9) yields
  \begin{equation}
  f''\cdot\frac{|\nabla u|^{2}+4u_{0}u_{i}\phi_{i}}{4u_{0}(1+|\nabla
  u|^{2})}+f'\cdot g^{ij}\phi_{ij}+\frac{\eta R_{kj}\nu^{k}\nu^{j}}{v^{2}}\leq0.
  \end{equation}

  {\bf Step 3.} Note that
  $$
  \phi_{ij}=(-\rho^{2})_{ij}=-2\rho_{i}\rho_{j}-2\rho\rho_{ij}.
  $$
  So
  $$
  g^{ij}\phi_{ij}=-2(\rho_{i}^{2}-\rho_{i}\nu^{i}\rho_{j}\nu^{j})-2\rho(\rho_{ii}-\rho_{ij}\nu^{i}\nu^{j})
  $$
  $$
  \hspace{11mm}=-2(1-\rho_{i}\nu^{i}\rho_{j}\nu^{j})-2\rho(\triangle\rho-\rho_{ij}\nu^{i}\nu^{j}).
  $$
  For
  $$
  |\rho_{i}\nu^{i}|\leq1,|\rho_{i}\nu^{i}\rho_{j}\nu^{j}|\leq1,
  $$
  we get
  \begin{equation}
  -2(1-\rho_{i}\nu^{i}\rho_{j}\nu^{j})\geq-2(1-(-1))=-4.
  \end{equation}

  Next we want to estimate $\rho_{ij}$. Denote the Hessian of
  $\rho$ by $H(\rho)$. Also denote $\frac{\partial}{\partial
  x^{i}}$ by $X_{i}$. Let $\widetilde{X}_{i}$, $i=1,\cdots,n$, be
  the Jacobi fields along $\gamma$ satisfying
  $$
  \widetilde{X}_{i}(\gamma(\rho))=X_{i}(\gamma(\rho)),
  \widetilde{X}_{i}(\gamma(0))=X_{i}(\gamma(0)),
  [\widetilde{X}_{i},\frac{\partial}{\partial\gamma}]=0.
  $$
  Since $\widetilde{X}_{i}$ is a Jacobi field, it satisfies the
  Jacobi equation
  $$
  \nabla_{\frac{\partial}{\partial\gamma}}\nabla_{\frac{\partial}{\partial\gamma}}
  \widetilde{X}_{i}+R(\widetilde{X}_{i},\frac{\partial}{\partial\gamma})\frac{\partial}{\partial\gamma}
  =0.
  $$
  We compute
  $$
  \hspace{-32mm}\rho_{ij}=H(\rho)(X_{i},X_{j})
  $$
  $$
  \hspace{-16mm}=X_{i}X_{j}\rho-(\nabla_{X_{i}}X_{j})\rho
  $$
  $$
  \hspace{-3mm}=X_{i}\langle X_{j},\frac{\partial}{\partial\rho}\rangle-\langle\nabla_{X_{i}}X_{j},
  \frac{\partial}{\partial\rho}\rangle
  $$
  $$
  \hspace{-27mm}=\langle X_{j},\nabla_{X_{i}}\frac{\partial}{\partial\rho}\rangle
  $$
  $$
  \hspace{-28mm}=\langle X_{j},\nabla_{\frac{\partial}{\partial\rho}}X_{i}\rangle
  $$
  $$
  \hspace{-15mm}=\int_{0}^{\rho}\frac{d}{dt}\langle\widetilde{X}_{j},\nabla_{\frac{\partial}{\partial t}}
  \widetilde{X}_{i}\rangle dt
  $$
  $$
  \hspace{22mm}=\int_{0}^{\rho}(\langle\nabla_{\frac{\partial}{\partial t}}
  \widetilde{X}_{j},\nabla_{\frac{\partial}{\partial t}}
  \widetilde{X}_{i}\rangle+\langle\widetilde{X}_{j},
  \nabla_{\frac{\partial}{\partial t}}\nabla_{\frac{\partial}{\partial t}}
  \widetilde{X}_{i}\rangle)dt
  $$
  $$
  \hspace{24mm}=\int_{0}^{\rho}(\langle\nabla_{\frac{\partial}{\partial t}}
  \widetilde{X}_{j},\nabla_{\frac{\partial}{\partial t}}
  \widetilde{X}_{i}\rangle-\langle\widetilde{X}_{j},
  R(\widetilde{X}_{i},\frac{\partial}{\partial t})\frac{\partial}{\partial t}\rangle)dt
  $$
  Note that $X_{n}=\frac{\partial}{\partial \rho}$. Then the fourth equality shows that
  $$
  \rho_{nj}=0
  $$
  for $j=1,\cdots,n$.

  Define the "index form" for a vector field $\widetilde{X}_{i}$
  along $\gamma$ (\cite{CE}) as follows
  $$
  I_{0}^{\rho}(\widetilde{X}_{i}):=\int_{0}^{\rho}(|\nabla_{\frac{\partial}{\partial t}}
  \widetilde{X}_{i}|^{2}-\langle\widetilde{X}_{i},
  R(\widetilde{X}_{i},\frac{\partial}{\partial t})\frac{\partial}{\partial
  t}\rangle)dt.
  $$
  By this definition,
  $$
  \rho_{ii}=I_{0}^{\rho}(\widetilde{X}_{i}).
  $$

  Let $E_{i}(\gamma(t))$, $0\leq t\leq\rho$, $i=1,\cdots,n$, be
  the parallel transport of $X_{i}(\gamma(\rho))$ along $\gamma$.
  Since Jacobi field minimizes the index form among all vector
  fields along the same geodesic with the same boundary values (see
  \cite{CE}), we have
  $$
  \hspace{-33mm}\rho_{ii}\leq I_{0}^{\rho}(\frac{t}{\rho}E_{i})
  $$
  $$
  =\int_{0}^{\rho}(\frac{1}{\rho^{2}}-\frac{t^{2}}{\rho^{2}}K(E_{i},\frac{\partial}{\partial
  t}))dt
  $$
  $$
  =\frac{1}{\rho}-\frac{1}{\rho^{2}}\int_{0}^{\rho}t^{2}K(E_{i},\frac{\partial}{\partial
  t})dt.
  $$
  It follows that
  $$
  \hspace{-48mm}\triangle\rho=\sum_{i=1}^{n-1}\rho_{ii}
  $$
  $$
  =\frac{n-1}{\rho}-\frac{1}{\rho^{2}}\int_{0}^{\rho}t^{2}Ric(\frac{\partial}{\partial
  t},\frac{\partial}{\partial
  t})dt
  $$
  $$
  \hspace{-5mm}\leq\frac{n-1}{\rho}+\frac{1}{\rho^{2}}\int_{0}^{\rho}t^{2}(n-1)Qdt
  $$
  $$
  \hspace{-18mm}=\frac{n-1}{\rho}+\frac{(n-1)Q}{3}\rho.
  $$

  On the other hand,
  $$
  -\rho_{ij}\nu^{i}\nu^{j}=-\int_{0}^{\rho}\langle\nabla_{\frac{\partial}{\partial t}}\widetilde{X}_{j},
  \nabla_{\frac{\partial}{\partial
  t}}\widetilde{X}_{i}\rangle\nu^{i}\nu^{j}dt+\int_{0}^{\rho}\langle\nu^{j}\widetilde{X}_{j},
  R(\nu^{i}\widetilde{X}_{i},\frac{\partial}{\partial t})\frac{\partial}{\partial
  t}\rangle dt.
  $$
  The first term on the r.h.s. can be bounded as follows
  $$
  \hspace{-32mm}-\int_{0}^{\rho}\langle\nabla_{\frac{\partial}{\partial t}}\widetilde{X}_{j},
  \nabla_{\frac{\partial}{\partial t}}\widetilde{X}_{i}\rangle\nu^{i}\nu^{j}dt
  $$
  $$
  \hspace{-45mm}\leq\int_{0}^{\rho}|\langle\nabla_{\frac{\partial}{\partial t}}\widetilde{X}_{j},
  \nabla_{\frac{\partial}{\partial t}}\widetilde{X}_{i}\rangle|dt
  $$
  $$
  \hspace{-46mm}\leq\int_{0}^{\rho}|\nabla_{\frac{\partial}{\partial t}}\widetilde{X}_{j}|\cdot
  |\nabla_{\frac{\partial}{\partial t}}\widetilde{X}_{i}|dt
  $$
  $$
  \hspace{-23mm}\leq(\int_{0}^{\rho}|\nabla_{\frac{\partial}{\partial t}}\widetilde{X}_{j}|^{2}dt)^{\frac{1}{2}}\cdot
  (\int_{0}^{\rho}|\nabla_{\frac{\partial}{\partial t}}\widetilde{X}_{i}|^{2}dt)^{\frac{1}{2}}
  $$
  $$
  \hspace{-26mm}\leq\frac{1}{2}\int_{0}^{\rho}|\nabla_{\frac{\partial}{\partial t}}\widetilde{X}_{j}|^{2}dt
  +\frac{1}{2}\int_{0}^{\rho}|\nabla_{\frac{\partial}{\partial t}}\widetilde{X}_{i}|^{2}dt
  $$
  $$
  \hspace{-49mm}=(n-1)\int_{0}^{\rho}|\nabla_{\frac{\partial}{\partial t}}\widetilde{X}_{j}|^{2}dt
  $$
  $$
  \hspace{-16mm}=(n-1)(I_{0}^{\rho}(\widetilde{X}_{j})+\int_{0}^{\rho}\langle\widetilde{X}_{j},R(\widetilde{X}_{j},
  \frac{\partial}{\partial t})\frac{\partial}{\partial t}\rangle dt)
  $$
  $$
  \hspace{-13mm}\leq(n-1)(I_{0}^{\rho}(\frac{t}{\rho}E_{j})+\int_{0}^{\rho}\langle\widetilde{X}_{j},R(\widetilde{X}_{j},
  \frac{\partial}{\partial t})\frac{\partial}{\partial t}\rangle dt)
  $$
  $$
  \hspace{10mm}\leq\frac{(n-1)^{2}}{\rho}+\frac{(n-1)^{2}Q}{3}\rho+(n-1)\int_{0}^{\rho}\langle\widetilde{X}_{j},R(\widetilde{X}_{j},
  \frac{\partial}{\partial t})\frac{\partial}{\partial t}\rangle
  dt.
  $$

  To go further, we need to estimate $|\widetilde{X}_{j}|$. This
  is done in the next step.

  {\bf Step 4.} Generally, let $J=J(t)$ be a Jacobi field along
  $\gamma$ satisfying
  $$
  |J(0)|=0,|J(\rho)|=1.
  $$
  We compute
  $$
  \frac{d}{dt}|J(t)|^{2}=2\langle J,\nabla_{\frac{\partial}{\partial
  t}}J\rangle,
  $$
  $$
  \frac{d^{2}}{dt^{2}}|J(t)|^{2}=2|\nabla_{\frac{\partial}{\partial
  t}}J|^{2}+2\langle J,\nabla_{\frac{\partial}{\partial
  t}}\nabla_{\frac{\partial}{\partial
  t}}J\rangle
  $$
  $$
  \geq-2\langle J,R(J,\frac{\partial}{\partial
  t})\frac{\partial}{\partial
  t}\rangle
  $$
  $$
  \hspace{-13mm}\geq-2|J|^{2}C_{+},
  $$
  i.e.,
  $$
  \frac{d^{2}}{dt^{2}}|J(t)|^{2}+2C_{+}|J(t)|^{2}\geq0.
  $$

  Define
  $$
  L:=\frac{d^{2}}{dt^{2}}+2C_{+}
  $$
  be an ordinary differential operator. Let
  $$
  w(t)=\cos\sqrt{2C_{+}}t.
  $$
  Then $w$ satisfies
  $$
  Lw=0.
  $$
  Moreover, note that $0\leq C_{+}<\frac{\pi^{2}}{8r^{2}}$, so $w$
  also satisfies
  $$
  w(t)>0,0\leq t\leq\rho.
  $$

  By Theorem 2.11 in \cite{L}, we have
  $$
  \max_{0\leq t\leq\rho}\{\frac{|J(t)|^{2}}{w(t)}\}=\max\{\frac{|J(0)|^{2}}{w(0)},\frac{|J(\rho)|^{2}}{w(\rho)}\}
  $$
  $$
  \hspace{11mm}=\max\{0,\frac{1}{w(\rho)}\}
  $$
  $$
  \hspace{-3mm}=\frac{1}{w(\rho)}.
  $$
  Hence
  $$
  |J(t)|^{2}\leq\frac{w(t)}{w(\rho)}=\frac{\cos\sqrt{2C_{+}}t}{\cos\sqrt{2C_{+}}\rho}
  <\frac{1}{\cos\sqrt{2C_{+}}r}.
  $$

  {\bf Step 5.} Applying the results in {\bf Step 4} to
  $\widetilde{X}_{i}$, we get
  $$
  -\rho_{ij}\nu^{i}\nu^{j}\leq\frac{(n-1)^{2}}{\rho}+\frac{(n-1)^{2}Q}{3}\rho
  +(n-1)\int_{0}^{\rho}|\widetilde{X}_{j}|^{2}K(\widetilde{X}_{j},\frac{\partial}{\partial
  t})dt
  $$
  $$
  \hspace{-24mm}+\int_{0}^{\rho}|\nu^{i}\widetilde{X}_{i}|^{2}K(\nu^{i}\widetilde{X}_{i},\frac{\partial}{\partial
  t})dt
  $$
  $$
  \hspace{-9mm}\leq\frac{(n-1)^{2}}{\rho}+\frac{(n-1)^{2}Q}{3}\rho+\frac{n(n-1)C_{+}}{\cos\sqrt{2C_{+}}r}\rho.
  $$
  Therefore
  $$
  \hspace{-76mm}\triangle\rho-\rho_{ij}\nu^{i}\nu^{j}
  $$
  $$
  \leq\frac{n-1}{\rho}+\frac{(n-1)Q}{3}\rho
  +\frac{(n-1)^{2}}{\rho}+\frac{(n-1)^{2}Q}{3}\rho+\frac{n(n-1)C_{+}}{\cos\sqrt{2C_{+}}r}\rho
  $$
  \begin{equation}
  \hspace{-35mm}=\frac{n(n-1)}{\rho}+\frac{n(n-1)Q}{3}\rho+\frac{n(n-1)C_{+}}{\cos\sqrt{2C_{+}}r}\rho.
  \end{equation}

  Combining (11) and (12) yields
  $$
  g^{ij}\phi_{ij}
  \geq-4-2n(n-1)-\frac{2n(n-1)Q}{3}\rho^{2}-\frac{2n(n-1)C_{+}}{\cos\sqrt{2C_{+}}r}\rho^{2}
  $$
  $$
  \hspace{9mm}\geq-4-2n(n-1)-\frac{2n(n-1)Q}{3}r^{2}-\frac{2n(n-1)C_{+}}{\cos\sqrt{2C_{+}}r}r^{2}
  $$
  \begin{equation}
  \hspace{-66mm}:=-C_{2}.
  \end{equation}

  {\bf Step 6.} When $|\nabla u|(P)\geq16u_{0}$, we have
  $$
  |\nabla u|^{2}+4u_{0}u_{i}\phi_{i}\geq|\nabla u|^{2}-8|\nabla u|u_{0}\geq\frac{1}{2}|\nabla u|^{2}.
  $$
  Moreover, when $|\nabla u|(P)\geq\max\{3,16u_{0}\}$, we have
  \begin{equation}
  \frac{|\nabla u|^{2}+4u_{0}u_{i}\phi_{i}}{4u_{0}^{2}(1+|\nabla u|^{2})}
  \geq\frac{1}{8u_{0}^{2}}\cdot\frac{|\nabla u|^{2}}{1+|\nabla u|^{2}}
  >\frac{1}{10u_{0}^{2}}.
  \end{equation}

  Note that
  \begin{equation}
  \frac{R_{kj}\nu^{k}\nu^{j}}{v^{2}}\geq-\frac{(n-1)Q|\nu|^{2}}{v^{2}}\geq-(n-1)Q,
  \end{equation}
  and
  $$
  f'=C_{1}\exp(C_{1}\phi),f''=C_{1}^{2}\exp(C_{1}\phi).
  $$
  Combing (10) and (13)-(15), we get
  \begin{equation}
  \frac{1}{10u_{0}^{2}}C_{1}^{2}-C_{2}C_{1}-(n-1)Q\leq0.
  \end{equation}
  It follows that for large $C_{1}$ (depending only on $n$, $Q$, $r$, $C_{+}$, $u_{0}$),
  (15) is contradicted if $|\nabla u|(P)\geq\max\{3,16u_{0}\}:=C_{3}$. Therefore
  $$
  |\nabla u|(P)\leq C_{3}.
  $$
  $$
  v(P)\leq 1+C_{3}:=C_{4}.
  $$
  $$
  \eta(x,u(x))v(x)\leq\eta(P,u(P))v(P)\leq C_{4}\exp(C_{1}r^{2}).
  $$
  At point $p$,
  $$
  (\exp(r^{2}-\frac{1}{2})-1)v(p)\leq C_{4}\exp(C_{1}r^{2}).
  $$
  For large $u_{0}$ it is easy to see that $C_{1}$ may be taken to be a multiple of $u_{0}^{2}$,
  so that the interior gradient bound has the form
  $$
  v(p)\leq K_{2}\exp(K_{1}u_{0}^{2}),
  $$
  where $K_{i}$, $i=1,2$, depend only on $n$, $Q$, $C_{+}$, $r$.

  This completes the proof of Theorem 2.
\end{section}


\begin{thebibliography}{10}
\bibitem{BDM}
E. Bombieri, E. De Giorgi and M. Miranda, {\em Una maggiorazione a
priori relativa alla ipersuperfici minimali non parametriche},
{\em Arch. Rational Mech. Anal.}, {\bf 32}(1969)255-267.

\bibitem{CE}
J. Cheeger and D. G. Ebin, {\em Comparison Theorems in Riemannian
Geometry}, North Holland, Amsterdam, 1975.

\bibitem{CM}
T. H. Colding and W. P. Minicozzi II, {\em Sharp estimates for
mean curvature flow of graphs}, {\em arXiv:math.AP/0305099 v2}.

\bibitem{EH} Ecker and Huisken, {\em Interior estimates for hypersurfaces
moving by mean curvature}. Invent. Math. 105 (1991), no. 3,
547--569.

\bibitem{K}
N. J. Korevaar, {\em An easy proof of the interior gradient bound
for solutions to the prescribed mean curvature equation}, {\em
Proc. Symp. Pure Math.}, Part !!, 1986,45:81-89.

\bibitem{L}
F. H. Lin, Q. Han, {\em Elliptic partial differential equations},
New York : Courant Institute of Mathematical Sciences ; 2000.

\bibitem{M1}
L. Ma, {\em On minimal graph evolutions in the hyperbolic space},
{\em Acta Math. Sinica}, English Series, vol. {\bf 15},
No.3(1999), 371-374.

\bibitem{M2}
L. Ma, {\em Mean curavture flow in a product manifold}. in
preparation.

\bibitem{R}
H. Rosenberg, {\em On minimal surface in $M\times R$}, Preprint,
2001.


\bibitem{SY}
R. Schoen and S.-T. Yau, {\em Lectures on Differential Geometry},
International Press Inc., Boston, 1994.

\end{thebibliography}
\end{document}